\documentclass[12pt]{article}
\usepackage{amssymb,latexsym}
\usepackage{amsmath}
\usepackage{amsfonts}

\usepackage{graphicx}
\usepackage[all,dvips]{xy}

\begin{document}

\title{Oppositions in a line segment}
      
\date{}

\author{\small Alexandre Costa-Leite\\
\small Department of Philosophy\\
\small University of Brasilia, Brazil}

\maketitle

\begin{abstract}
Traditional oppositions are at least two-dimensional in the sense
that they are built based on a famous bidimensional object called
\emph{square of oppositions} and on one of its extensions such as
Blanch\'e's hexagon.
Instead of two-dimensional objects, 
this article proposes a construction to deal with oppositions
in a one-dimensional line segment. 
\end{abstract}

%%%%%%%%%%%% Content of the article %%%%%%%%%%%%

\section*{Introduction}

                The basic theory of oppositions has been developed
considering a bidimensional structure, i.e., the well known
\emph{square of oppositions}. Later, it has been generalized
to an hexagon of oppositions by Blanch\'e in \cite{blanche}.
Moretti argues in \cite{moretti1} that there is a `geometry' of oppositions,
generalizing squares and hexagons, basically, to three-dimensional
structures such as cubes and tetradecahedrons. This
gives rise to the domain of \emph{n}-opposition theory.
Since works proposed by Beziau in \cite{beziau2} and Moretti 
in \cite{moretti1}, there are now many
researches in the field . The reader should
check \cite{parsons} for an introduction to the square
of oppositions and \cite{moretti2} to its main recent developments. Beziau
and Read stated in \cite{beziau-read} that the theory of oppositions
cannot be identified with the diagram representing this theory (p.315).
This means that there are much more on oppositions than what is
represented in the relations between corners of the square or the hexagon.
This article\footnote{A previous version of this paper appeared as
a preprint in arXiv (April, 2016). Thanks to Rodrigo Freire, Edelcio de Souza and Fabien
Schang for remarks on the constructions proposed here.}, in some sense, can be viewed as an attempt to justify
the claim that the theory of oppositions is not only the study of these
\emph{n}-dimensional diagrams ($n \geq 2$).

Consider a question: is there a way to represent oppositions without 
two-dimensional objects such as squares or objects of
higher dimensions? The answer is \emph{yes}. A construction
to formulate this reduction is proposed showing that there is no
need for two-dimensional objects to establish the basic
theory of oppositions. Indeed, one dimension is enough.
This means that oppositions can be defined
in a line segment, a piece of one-dimensional space. Moreover, oppositions
require precisely at least one dimension to be defined, 
and the traditional case of the bidimensional square can be converted step-by-step to it. 
However, applying the same strategy to Blanch\'e's hexagon does not work. But
the situation changes with some constraints added at the level of the  basic construction. So, it is also possible to convert the standard bidimensional Blanch\'e's hexagon to a line segment.

This paper shows that line segments are
sufficient to define four basic standard oppositions in such a way that the
square can be derived from this primitive structure. This line segment of oppositions 
is called here \emph{basic construction}. Then, this
very same basic construction does not work when applied to the standard hexagon (i.e. Blanch\'e's hexagon),  but there is a way to generate a similar strategy \emph{mutatis
mutandis} to it.  In what follows these constructions are explained 
presenting their range and limits.

\section{Defining oppositions in a line segment}

\subsection{The square}

Assume classical logic.
Given the framework of first-order logic, the square of oppositions uses four kinds
of categorical propositions (where $\varphi$ is formula): (A) universal affirmative of the form
$\forall x \varphi$, (E) universal negative of the form
$\forall x \neg \varphi$, (I) existential affirmative $\exists x \varphi$
and (O) existential negative $\exists x \neg \varphi$. 
Taking into account that $\forall$ and $\exists$
are interdefinable in the presence of negation, it is a matter of taste to decide
which one to use to represent these four propositions. They
can appear as above, mixing both quantifiers, or with only
one kind of quantifier. Thus, for universal quantifier
and negation there is the following: (A) $\forall x \varphi$, (E) $\forall x \neg \varphi$,
(I) $\neg \forall x \neg \varphi$ and (O) $\neg \forall x \varphi$.
For existential quantification
and negation: (A) $\neg \exists x \neg \varphi$, (E) $\neg \exists x \varphi$,
(I) $\exists x \varphi$ and (O) $\exists x \neg \varphi$. In general, in the literature,
there are four traditional oppositions holding between these propositions: contradiction (d),
contrariety (c), subcontrariety (sc) and subalternation (s), which 
are defined in the standard way (see \cite{parsons}).

Further, there are, indeed, other families
of concepts satisfying oppositional structures:
they could be metaphysical statements
such as \emph{necessity} and \emph{possibility} 
(and their derivatives), or deontic propositions containing notions of \emph{obligation}
and \emph{permission}, or even statements containing temporal aspects 
such as \emph{always}, \emph{sometimes} and \emph{never}. So, a pure oppositional structure is not necessarily decorated with categorical statements. Indeed, there are many possible decorations
of the square (see \cite{blanche, beziau, moretti1}). These concepts fit pretty
well satisfying the structure of the square of oppositions (or its extensions). 
They - and similar concepts - are here called categorical-like concepts: these
are notions satisfying the four oppositions and, for this reason, can be arranged inside
the traditional square.  A set of categorical-like
statements is denoted by $\mathbf{C}$.  
Oppositions are usually (and historically) represented bidimensionally
using a two-dimensional object with
\emph{lenght} and \emph{width} (i.e. a square):

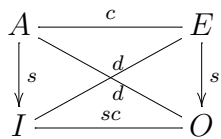
\[  \xymatrix{
  A \ar@{-}[rr]^{c} \ar@{->}[d]^{s} \ar@{-}[rrd]^{d} && E \ar@{->}[d]^{s} \ar@{-}[dll]^{d} \\
    I \ar@{-}[rr]^{sc} && O
   }\] \\

I argue that the above two-dimensional square can be reduced to a simple one-dimensional
structure, that is, a line segment of integers.\footnote{\emph{Simple} in the sense that it requires
only one dimension.} To show how to build
this reduction is the \emph{first construction} suggested here.

Assume the set $\mathbb{Z}$ of integers and, then, the sets $\mathbb{Z_{+}}$ (positive integers) and
$\mathbb{Z_{-}}$ (negative integers). If the set of integers does not have
zero, it is denoted by $\mathbb{Z^{*}}$, as usual. 
Take a line segment such that for each $j \in \mathbb{Z}$, $j \neq 0$, there exists $-j$, the
symmetric of $j$. In particular, consider a set $\mathbb{Z'}= \{-r,-q,q,r\} \subseteq \mathbb{Z}$.
Then the procedure is: a categorical-like proposition 
is assigned to each element $j \in \mathbb{Z'}$ in the following way:

\begin{itemize}
\item $j \in \mathbb{Z^{*}_{+}}$ if, and only if, $j$ is associated to a universal statement;

\item $j \in \mathbb{Z^{*}_{-}}$ if, and only if, $j$ is associated to an existential statement.

\end{itemize}

There is a division of propositions in the sets of universal and existential statements. Let $\mathbf{C}$ be a
set of categorical (or categorical-like) propositions. The function $i$ which
connects elements of $\mathbf{C}$ to elements of $\mathbb{Z'}$ is a bijection such that: $i(A) = q$;
$i(E) = r$; $i(I) = -r$ and $i(O) = -q$.

We use $\alpha, \beta$ for arbitrary propositions ranging on $\mathbf{C}$. Traditional oppositions
now have to be reformulated inside this new framework. So, clauses for oppositions are defined as follows:

\begin{enumerate}
\item Propositions $\alpha, \beta$ are \emph{contradictories} if, and only if,
their assigned numbers have sum equals to $0$, i.e, $i(\alpha) + i(\beta)=0$;

\item Propositions $\alpha, \beta$ are \emph{contraries} if, and only if,
$i(\alpha), i(\beta) \in \mathbb{Z^{*}_{+}}$;

\item Propositions $\alpha, \beta$ are \emph{subcontraries} if, and only if,
$i(\alpha), i(\beta) \in \mathbb{Z^{*}_{-}}$;

\item $\beta$ is subaltern of $\alpha$ if, and only if, $i(\beta) \neq -i(\alpha)$ and $i(\beta) \in \mathbb{Z^{*}_{-}}$;

\end{enumerate}

This completes the construction. Let's show with an example how the square can be defined in this way: consider,
for instance, that $\mathbb{Z'}= \{-2,-1,\allowbreak 1,2\} \subseteq \mathbb{Z}$. Thus, given that $1, 2 \in \mathbb{Z^{*}_{+}}$, 
it follows that these numbers are assigned to universal propositions (in which way this association
is done is not important). In
the same manner, from the fact that $-2, -1 \in \mathbb{Z^{*}_{-}}$, these numbers are connected
to existential statements. Without loss generality, 
consider that $1$ is the number connected to (A), i.e, $i(A)=1$ and $i(E)=2$.
Therefore, the contradictory of (A) is the proposition (O) because $i(O)=-1$ satisfying, therefore, the condition
to be a contradiction, i.e., $1 + (-1)=0$. The same for the relations between
(E) and (I).
By construction, (A) and (E) are contraries while (I) and (O) are subcontraries, and these
are subalterns of (A) and (E), respectively. (I) is subaltern of (A) given that (I) has
associated to it the integer $-2$ which is, in its turn, different of $-1$ (note that $-i(A)=i(O)$).
As far as there are only four oppositions, this procedure can always be done, no matter
the family of concepts considered.
This is an example provided to show how
a two-dimensional square of oppositions can be reduced to a one-dimensional line segment of oppositions.

Despite the beauty of oppositions represented in a square, these oppositions can be converted into a one-dimensional line segment,
if the number of oppositions remain four. This is the main contribution of this paper. 
Nevertheless, not all two-dimensional objects can be reduced using
this same strategy. It already fails in the case of the hexagon. So, the basic
construction has to be improved in order to work also for the hexagon. This is of secondary
importance here, given that the hexagon does not have the same historical relevance of the square.

\subsection{The hexagon}

Blanch\'e proposed an extension of the square, the hexagon, and this new tool (still two-dimensional) can be used to model many situations
(see \cite{blanche}, but also \cite{beziau}). 
Following Blanch\'e's construction, in the (incomplete) diagram below, (U) is defined as the disjunction $A \vee E$ and (Y) is the conjunction $I \wedge O$. 

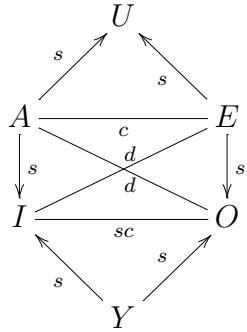
 \[ \xymatrix{
    & U \\
    A \ar@{->}[ur]^{s} \ar@{->}[d]^{s} \ar@{-}[rrd]^{d} && E \ar@{->}[ul]^{s} \ar@{-}[ll]^{c} \ar@{-}[lld]^{d}  \ar@{->}[d]^{s}  \\
    I && O \ar@{-}[ll]^{sc}\\
 &Y \ar@{->}[ul]^{s} \ar@{->}[ur]^{s} 
}\] \\

How can the reduction strategy, i.e., the basic construction be also extended 
to Blanch\'e's hexagon?
The fact that there are propositions of the form (U) and (Y) require some adaptations in the basic construction.
The \emph{second construction} consists in executing these adaptations in order to transform also the hexagon into a one-dimensional structure.

Take a line segment as above and 
a set $\mathbb{Z''}= \{-s,-r,-q,q,r,s\} \subseteq \mathbb{Z}$.
Now, as there are also a disjunction and a conjunction, some changes have
to be made in the way integers are assigned to propositions. To each element $j \in \mathbb{Z''}$, 
a proposition is connected as follows:

\begin{itemize}
\item  $j \in \mathbb{Z^{*}_{+}}$ if, and only if, $j$ is associated to a universal statement or a disjunction;

\item  $j \in \mathbb{Z^{*}_{-}}$ if, and only if, $j$ is associated to an existential statement or a conjunction;.

\end{itemize}

The division of the class of propositions now contains a set of universal and disjunctive statements, 
from one side, and a set of existential and conjunctive statements, from other side. Let the set $\mathbf{C'}$ be an 
expansion of the set $\mathbf{C}$ defined by
categorical-like propositions plus a disjunction and a conjunction, and let the function $i$ which
connects elements of $\mathbf{C'}$ to elements of $\mathbb{Z''}$ be a bijection such that: $i(U) = s$,
$i(A) = q$, $i(E) = r$, $i(I) = -r$, $i(O) = -q$ and $i(Y) = -s$. The number assigned to the conjunction of (I) and (O), that is, (Y), as well
to the disjunction of (A) or (E), that is (U), is the sum of both conjuncts (in the first case)
and the sum of both disjuncts (in the second case). So, define that the integers associated to statements (U) and (Y) are obtained
by the sums of their components: $i(U) = i(A) + i(E)$ and $i(Y) = i(I) + i(O)$. These numbers obtained by sums play an important role and, therefore, they are called
\emph{distinct objects}: the first one is the positive distinct object and the second is the negative. 
Consider a
similar strategy as the basic construction and let's try to reduce the hexagon to a line segment. 
Suppose,
for instance, that $\mathbb{Z''}= \{-3,-2,-1,1,2,3\} \subseteq \mathbb{Z}$ is given and numbers
are associated to propositions: $i(A)=1$, $i(E)=2$, $i(O)=-1$ and $i(I)=-2$. Thus, distinct
objects have the following association: $i(U)=3$ and $i(Y)=-3$, if we consider this particular
line segment $[-3, 3]$. Hence, these statements are
obviously in the opposition of contradiction. In this sense, clause (1) above, i.e., contradiction, holds for (U) and (Y):\footnote{Moretti (see \cite{moretti2}) remarked that the opposition
of contradiction can be characterized, no matter which dimension is considered, as a certain kind of symmetry. It is a conjecture of this work that 
the whole of \emph{n}-opposition theory can be reduced to variations of the constructions presented
in this paper.}

  \[ \xymatrix{
 U \\
  Y \ar@{-}[u]^{d} \\
}\] \\

But not all clauses from 1-4 are valid. 
The reason for this is that Blanch\'e's hexagon contains some unexpected connections
between propositions (Y)-(A)-(E) (they are contraries and then clause
2 fails):

 
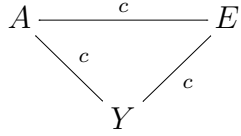
\[\xymatrix{
    A \ar@{-}[rr]^{c} \ar@{-}[rd]^{c} && E \ar@{-}[ld]^{c} \\ & Y 
   }\] \\

and propositions (U)-(I)-(O) (they are subcontraries and
then clause 3 fails):

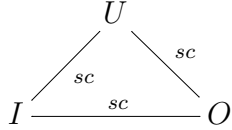
\begin{displaymath}
    \xymatrix{ & U \ar@{-}[rd]^{sc} \ar@{-}[ld]^{sc}\\
                I \ar@{-}[rr]^{sc} && O  }
\end{displaymath}

In addition, subalternation also fails. Thus, it is not straighforward to settle Blanch\'e's hexagon in a line segment of integers.\footnote{Note that these last three diagrams are displayed in Blanch\'e's hexagon (see \cite{blanche}). When these relations are added to the diagram, we have the complete
hexagon of oppositions.} So, 
some adaptations and repairs in the basic construction should have to be done.

Assume that $\alpha, \beta, \gamma$ are letters for arbitrary propositions. 
The construction below shows how to reduce Blanch\'e's hexagon to a line segment.
Let $i(\gamma)$ be a distinct object.
For contradictories, 
condition 1 remains the same as above, but for other oppositions, 
clauses are the following (note that they have to be applied
in this order: contradiction, contraries or subcontraries, and, by the end, 
subalternation. Then gaps will be gradually filled, as the
oppositions are excludent, i.e, two propositions cannot be related
by two different oppositions):

\begin{enumerate}

\item[1*.] Propositions $\alpha,\beta, \gamma$ are \emph{contraries} if, and only if,
$(i(\alpha) + i(\beta)) + i(\gamma) =0$ and $i(\gamma) \in \mathbb{Z^{*}_{-}}$;

\item[2*.] Propositions $\alpha,\beta, \gamma$ are \emph{subcontraries} if, and only if,
$(i(\alpha) + i(\beta)) + i(\gamma) =0$ and $i(\gamma) \in \mathbb{Z^{*}_{+}}$;

\item[3*.] $\beta$ is subaltern of $\alpha$ if, and only if, $i(\beta) \neq -i(\alpha)$ and $i(\beta) \in \mathbb{Z^{*}_{-}}$,
or $i(\beta) \neq -i(\alpha)$ and a) $i(\beta)>i(\alpha)$ and $i(\alpha),i(\beta) \in \mathbb{Z^{*}_{+}}$ or b)  
$i(\beta)>i(\alpha)$ and $i(\alpha),i(\beta) \in \mathbb{Z^{*}_{-}}$ ($i(\beta)$, in the last condition,
is a distinct object);

\end{enumerate}

It is not difficult to provide an example of integers to show that the hexagon can be reduced to
it, as done in the case of the square. Assume first - without loss of generality - that $i(A)=1$.
So, $i(E)=2$ and, thus, the distinct object $i(U)=i(A)+i(E)=3$. Second, consider that $i(O)=-1$. It follows that $i(I)=-2$
and that the other distinct object $i(Y)=i(O)+i(I)=-3$. Consequently, (U) and (Y) are contradictories. Note that
$(i(O)+i(I))+i(U)=0$, so these are subcontraries, as $i(U) \in \mathbb{Z^{*}_{+}}$, i.e, it is the positive distinct
object. Moreover,
$(i(A)+i(E))+i(Y)=0$, so these are contraries, as $i(Y) \in \mathbb{Z^{*}_{-}}$, i.e, it is the negative distinct
object. Concerning subalternation, the square is contained in the hexagon, so (I) is subaltern of (A) and
(O) is subaltern of (E), as in the case of the square. For other relations of subalternation, note that (U) is
subaltern of (A) and (E), because $i(U) > i(A)$ and $i(U) > i(E)$  while $i(U), i(A), i(E) \in \mathbb{Z^{*}_{+}}$.
Differently, (I) and (O) are subalterns of (Y), given that $i(Y) < i(I)$ and $i(Y) < i(O)$, and $i(Y), i(I), i(O) \in \mathbb{Z^{*}_{-}}$.

This is one way to reduce the hexagon to a one-dimensional object, but
it is not argued that there is no other way. 
Note also that only Blanch\'e's hexagon is taken into consideration, 
because it is the traditional hexagon, although
there are some other available in the literature.

\section{Conclusion}

The standard square of opposition is a two-dimensional object
used to organize and manage categorical-like concepts. 
There are many notions which fit in the square 
and its regular extensions such as the 
two-dimensional hexagon, and its three-dimensional
correlates.

This paper proposed a construction which shows that
there is no need to develop two-dimensional squares to represent categorical-like
statements or concepts, instead a one-dimensional line segment based
on integers is enough. Although the success in the case of
the square, we have showed that the same technique does not immediately
work for Blanch\'e's hexagon and, therefore, it does require some
adaptations and expansions to transform also the hexagon into a one-dimensional object. 
While both squares and hexagons have formulations in line
segments of integers, these formulations are not straightforward
and intuitive, and thus it seems these results are not so 
effective and simple
as the pictorial effect of two-dimensional diagrams. However, from the theoretical
viewpoint, it is important to know that there are simple one-dimensional objects
able to accommodate theory of oppositions, although
the fact that maybe, at this level, they are not so manageable as squares and hexagons.
This is a clue that oppositions do not match with the research on diagrams, and
in this sense, what authors supported in \cite{beziau-read} seems
to be rather right:

\begin{quotation}
``The diagram has been very important in promoting the theory but the theory does not reduce
to the diagram. The theory started many centuries before the basic diagram was drawn and
developed beyond this diagram...The strength of this theory is that it is at the same time fairly 
simple but quite rich; it can be applied to many different kinds of proposition, and also to objects and concepts.
It can also be generalized in various manners, in particular, by constructing many 
different geometrical objects.'' (pgs. 315-316, in \cite{beziau-read})\\
\end{quotation}

Oppositions are relations between propositions
and cannot be identified, therefore, with the study of diagrams, though
these are useful tools to explain what oppositions are
and, in this way, can be largely applied especially for learning purposes.

In general, authors working on the square of opposition generally
start with two-dimensional constructions like squares and hexagons,
and then they jump to three-dimensional solids and so on.
It seems a novelty the construction proposed in this paper because
it shows that one does not need to go \emph{n}-dimensional (for $n \geq 2$)
to develop the theory of oppositions. 

There are, notwithstanding, some problems which remain open: the question to 
determine whether the same procedure can also be 
applied to solids and higher dimensions, as well as to more than four oppositions,
are very complicated and still have to investigated in detail in the domain of 
line segments of oppositions.

%%%%%%%%%%%% Bibliography %%%%%%%%%%%%

\end{document}